\documentstyle[12pt,twoside]{article}
\newtheorem{theorem}{\bf Theorem}[section]
\newtheorem{proposition}[theorem]{\bf Proposition}
\newtheorem{lemma}[theorem]{\bf Lemma}
\newtheorem{corollary}[theorem]{\bf Corollary}
\newtheorem{example}[theorem]{\bf Example}

\newtheorem{remark}[theorem]{\bf Remark}

\input{amssym}

\date{}
\begin{document}

\title{{\Large\bf Cohomological properties of vector-valued Lipschitz algebras and their second duals}}

\author{{\normalsize\sc M. J. Mehdipour and A. Rejali\footnote{Corresponding author}}}
\maketitle

{\footnotesize  {\bf Abstract.} Let $\frak{F}(X, A)$ be one of the Banach algebras $\hbox{Lip}(X, A)$ or $\hbox{lip}(X, A)$.
In this paper, we show that $\frak{F}(X, A)$ is amenable if and only if $X$ is uniformly discrete and $A$ is amenable. We also prove that the result holds for $\hbox{lip}^\circ(X, A)$ instead of $\frak{F}(X, A)$. In the case where $A^*$ is separable, we establish that
$\frak{F}(X, A)^{**}$ is amenable if and only if $X$ is uniformly discrete and $A^{**}$ is amenable, however, amenability of $\hbox{lip}^\circ(X, A)^{**}$ is equivalent to amenability of $A^{**}$ and finiteness of $X$.
We prove that if $\hbox{Lip}(X, A)$ is point (respectively, weakly) amenable, then 
 $X$ is uniformly discrete and $A$ is point (respectively, weakly) amenable. In particular, $\hbox{Lip}X$ is weakly amenable if and only if $X$ is discrete. 
We  then investigate cohomological properties for vector-valued Banach algebras $C_0(X, A)$ and $L^1(G, A)$. Finally, we prove that biprojectivity (respectively, cyclically weak amenability) of $A^{**}$ implies biprojectivity (respectively, cyclically weak amenability)  of $A$. This result holds for weak amenability and cyclic amenability when $A$ is commutative.}
                     %-----------------------------------------
                     %-----------------------------------------
                     %-----------------------------------------
{}
                     %-----------------------------------------
                     %-----------------------------------------
                     %-----------------------------------------
  %                   \tableofcontents

\section{\normalsize\bf Introduction}

Let $X$ be a topological space and $(A, \|.\|)$ be a Banach algebra. A 
vector-valued function $f: X\rightarrow A$ is said to be \emph{bounded} if 
$$
\|f\|_{\infty, A}:=\sup_{x\in X}\|f(x)\|<\infty.
$$
We denote by $\ell^\infty(X, A)$ (respectively, $C(X, A)$) the Banach space consisting of all bounded (respectively, continuous) vector-valued functions from $X$ into $A$. Let $$C_b(X, A)= C(X, A)\cap \ell^\infty(X, A)$$ and $C_0(X, A)$ be the closed subalgebra of $C_b(X, A)$
consisting of all functions that vanish
at infinity.  It is known that  
$$
C_0(X,A)=
C_0(X)\check{\otimes}A,
$$
where ``$\check{\otimes}$" denotes the injective tensor product.

Now, let $(X, d)$ be a metric space, $f\in \ell^\infty(X, A)$ and $\alpha$ be a positive number. Then $f$ is called a \emph{vector-valued Lipschitz function of order} $\alpha$ if
$$
p_{d, A}(f):=\sup_{x, y\in X, x\neq y}\frac{\|f(x)-f(y)\|}{d^\alpha(x, y)}<\infty.
$$
The space of all vector-valued Lipschitz functions of order $\alpha$ from $(X, d)$ into $(A, \|.\|)$ is denoted by $\hbox{Lip}_\alpha(X, A)$.  It is proved that $\hbox{Lip}_\alpha(X, A)$ with the pointwise multiplication and the norm
$$
\|f\|_{d, A}:=\|f\|_{\infty, A}+p_{d, A}(f)
$$
is a Banach algebra.
Let $\hbox{lip}_\alpha(X, A)$ be the subalgebra of $\hbox{Lip}_\alpha(X, A)$ consisting of $f$ such that
$\lim_{r\rightarrow0}C_r=0$, where
$$
C_r=\sup\{\frac{\|f(x)-f(y)\|}{d^\alpha(x, y)}: x, y\in X, 0<d(x,y)<r\}.
$$
It is easy to see that
$$
\hbox{lip}_\alpha(X, A)\subseteq \hbox{Lip}_\alpha(X, A)\subseteq C_b(X, A)\subseteq\ell^\infty(X, A).
$$
If $X$ is a locally compact metric space, then we set $$\hbox{lip}^\circ_\alpha(X, A):=\hbox{lip}_\alpha(X, A)\cap C_0(X, A).$$
We write, $\hbox{Lip}(X, A)$, $\hbox{lip}(X, A)$ and $\hbox{lip}^\circ(X, A)$ for $\hbox{Lip}_1(X, A)$, $\hbox{lip}_1(X, A)$ and $\hbox{lip}^\circ_1(X, A)$, respectively. In the case where, $A=\Bbb{C}$, we write $\hbox{Lip}_\alpha(X, \Bbb{C})=\hbox{Lip}_\alpha X$, $\hbox{lip}_\alpha(X, \Bbb{C})=\hbox{lip}_\alpha X$ and $\hbox{lip}^\circ_\alpha (X, \Bbb{C})=\hbox{lip}^\circ_\alpha X$. 

Let us recall that $X$ is called \emph{uniformly discrete} if there exists $r>0$ such that $d(x, y)>r$ for all pairwise disjoint elements $x, y\in X$. One can prove that $X$ is uniformly discrete if and only if 
with equivalent norms we have
\begin{eqnarray*}
\hbox{lip}(X, A)=\ell^\infty(X,A)=
C(\beta X,A);
\end{eqnarray*}
in particular, if $X$ is finite, then there exists $n\in\Bbb{N}$ such taht with equivalent norms
\begin{eqnarray*}
\hbox{lip}(X,A)=
C(X,A)=
C(X)\check{\otimes}A=\sum_{i=1}^n\Bbb{C}\hat{\otimes} A,
\end{eqnarray*}
where $\beta X$, $\hat{\otimes}$ and $\|.\|_\epsilon$ denote the Stone-Cech compactification of $X$, the projective tensor product and the injective norm, respectively. Also, if $X$ is a locally compact metric space, then $X$ is uniformly discrete if and only if $\hbox{lip}^\circ(X, A)=c_0(X, A)$. For study of vector-valued Lipschitz algebras see \cite{abr, biyabani1, br}.

Let $A$ be a Banach algebra and $E$ be a Banach $A-$bimodule. Then a linear operator $D: A\rightarrow E$ is called a \emph{derivation} if $D(ab)=D(a)b+aD(b)$ for all $a, b\in A$. A derivation $D: A\rightarrow E$ is said to be \emph{inner} if there exists $z\in E$ such that $D(a)=az-za$ for all $a\in A$. Let $\varphi\in\Delta(A)$, the character space of $A$, and $E=\Bbb{C}$ with the following module action.
$$
\alpha\cdot a=a\cdot\alpha=\varphi(a)\alpha\quad(a\in A, \alpha\in\Bbb{C}).
$$
Then a continuous derivation $D: A\rightarrow \Bbb{C}$ is called a \emph{continuous point derivation at} $\varphi$ of $A$. Let us recall that $A$ is called \emph{amenable} if every continuous derivation from $A$ into $E^*$  is inner. If the only continuous derivation from $A$ into $A^*$ is inner derivations, then $A$ is called \emph{weakly amenable}. Also, $A$ is \emph{cyclically weakly amenable} if every continuous derivation $D: A\rightarrow A^*$ is cyclic,i.e., $\langle D(a), a\rangle=0$ for all $a\in A$. Finally, $A$ is \emph{point amenable} (respectively, 0-\emph{point amenable})  if there is no non-zero continuous point derivation at $\Delta(A)$ (respectively, $\Delta(A)\cup\{0\}$) of $A$. Clearly, every amenable Banach algebra is weakly amenable. It is proved that weak amenability implies cyclically weak amenability; also, cyclically weak amenability of $A$ force $A$ to be point amenable. For unital commutative Banach algebras, the concepts weak amenability, cyclically weak amenability and point amenability are equivalent. It is obvious that every point amenable Banach algebra is 0-point amenable. Also, if $A$ is essential, then point amenability and 0-point amenability are equivalent. One can prove that if $$\Delta(A)\neq\emptyset,$$ then $A$ is 0-point amenable if and only if $A$ cyclically weakly amenable; for details see \cite{mr4}.

Myers \cite{my} investigated the space of Lipschitz functions as a Banach algebra. But Sherbert's paper \cite{s2} is a fundamental contribution in the Lipschitz Banach algebras. He extended some results of Myers and studied ideals and point derivations in this Banach algebra. In the case where, $X$ is a non-discrete compact metric space and $\alpha\in (0, 1]$, Sherbert \cite{s2} showed that $\hbox{Lip}_\alpha X$ is not point amenable and thus it is neither weakly amenable nor cyclically weakly amenable. Bade, Curtis and Dales \cite{bcd} studied weak amenability of Lipschitz algebras when the metric space is compact. 
For infinite compact metric space $(X, d)$ and $\alpha\in (0, 1)$, they proved that $\hbox{lip}_\alpha X$ is always point amenable, however, if $\alpha\in (0, 1/2)$, then $\hbox{lip}_\alpha X$ is weakly amenable and so it is cyclically weakly amenable. Gourdeau \cite{g2} established that for any metric space $(X, d)$ and $\alpha>0$, the amenability of $\hbox{Lip}_\alpha X$, $\hbox{lip}_\alpha X$ or $\hbox{lip}_\alpha^\circ X$ is equivalent to uniform
discreteness of $X$; see also \cite{biyabani, cha, joh2, kaw,   s3, s4}.  Several authors were interested in the structure of space of vector-valued Lipschitz functions \cite{baca, gpr, iz, joh1, kxb}. In this paper, we continue these investigations on cohomological properties of vector-valued Lipschitz algebras.

In the following, let $X$ be a metric space, $Y$ is a locally compact metric space and $A$ is a Banach algebra. This paper is organized as follow.  In Section 2, we prove that $\hbox{Lip}(X, A)$ is amenable if and only if $\hbox{lip} (X, A)$ is amenable; or equivalent, $X$ is uniformly discrete and $A$ is amenable. We also show that $\hbox{lip}^\circ (Y, A)$ is amenable if and only if $Y$ is finite and $A$ is amenable. In the case where $A^*$ is separable, we show that 
$\hbox{Lip}(X, A)^{**}$ is amenable if and only if $X$ is uniformly discrete and $A^{**}$ is amenable, however, amenability of $\hbox{lip}^\circ(X, A)^{**}$ is equivalent to amenability of $A^{**}$ and finiteness of $X$.
 In Section 3, we prove that if $\hbox{Lip} (X, A)$ is point amenable, then $X$ is uniformly discrete and $A$ is point amenable. In the case where, $X$ is a discrete metric space and $A$ is point amenable, then every continuous point derivation at $X\otimes\Delta(A)$ of $\hbox{Lip} (X, A)$ is zero. In Section 4, we study weak amenability of $\hbox{Lip} X$ and show that $\hbox{Lip} X$ is weak amenable if and only if $X$ is discrete. Section 5 is devoted to study of point amenability of Banach algebras $C_0(X, A)$ and $L^1(G, A)$. We establish $C_0(X, A)$ is (respectively, 0-point) amenable if and only if $A$ is (respectively, 0-point) amenable, however, $L^1(G, A)$ is point amenable if and only if $A$ is point amenable. When $A$ is commutative, we establish that weak amenability of $C_0(X, A)$ and $A$ is equivalent. This result holds for $L^1(G, A)$.

\section{\normalsize\bf Amenability of Lipschitz vector-valued functions algebras}

We commence this section with the following lemma.

\begin{lemma}\label{havig1} Let $X$ and $Y$ be Banach spaces. Then $X\hat{\otimes}Y$ is a dense subset of $X\check{\otimes}Y$.
\end{lemma}
{\it Proof.} Since $X\otimes Y$ is a subset of $X\hat{\otimes}Y$ and $\|.\|_\epsilon\leq\|.\|_\pi$, it follows that $X\hat{\otimes}Y$ is a subset of $X\check{\otimes}Y$, where $\|.\|_\pi$ is the projective norm. Hence
$$
X\check{\otimes}Y=\overline{X\otimes Y}^{\|.\|_\epsilon}\subseteq\overline{X\hat{\otimes} Y}^{\|.\|_\epsilon}\subseteq\overline{X\check{\otimes} Y}^{\|.\|_\epsilon}=X\check{\otimes}Y.
$$
This shows that $X\check{\otimes}Y=\overline{X\hat{\otimes} Y}^{\|.\|_\epsilon}$.
$\hfill\square$\\

Gourdeau \cite{g1, g2} investigated amenability of Lipschitz vector-valued functions algebras. He \cite{g2} verified that if $(X, d)$ is any metric space, $\alpha$ is a positive real number and $A$ is commutative, then amenability of $\hbox{Lip}_\alpha (X, A)$, $\hbox{lip}_\alpha(X, A)$ or $\hbox{lip}_\alpha^\circ(X, A)$ implies uniform
discreteness of $X$. He in \cite{g1} showed that this result remains valid if $A$ be any Banach algebra and $\hbox{Lip}_\alpha (X, A)$, $\hbox{lip}_\alpha (X, A)$ or $\hbox{lip}_\alpha^\circ (X, A)$ separates the points of $X$. 

\begin{theorem}\label{poi}  Let $(X, d)$ be a metric space and $A$ be a Banach algebra with $\Delta(A)\neq\emptyset$.  Then the following assertions are equivalent.

\emph{(a)} $\emph{Lip}(X,A)$ is amenable.

\emph{(b)} $\emph{lip}(X,A)$ is amenable.

\emph{(c)} $X$ is uniformly discrete and $A$ is amenable.
\end{theorem}
{\it Proof.} Let $\frak{F}(X)$ be one of the Banach algebras $\hbox{Lip}X$ or $\hbox{lip}X$. Let also $\frak{F}(X, A)$ be the vector-valued Banach algebras corresponding to $\frak{F}(X)$. If  $x\in X$ and $\varphi\in\Delta(A)$, then the functions $f\mapsto x\otimes f$ and $f\mapsto \varphi\otimes f$ are continuous epimorphisms from $\frak{F}(X, A)$ onto $A$ and $\frak{F}(X)$, respectively, where 
$$
x\otimes f=f(x)\quad\hbox{and}\quad\varphi\otimes f(x)=\varphi(f(x));
$$
see Lemma 3.2 in \cite{mr6}.
So if $\frak{F}(X, A)$ is amenable, then the Banach algebras $\frak{F}(X)$ and $A$ are amenable. Hence $X$ is uniformly discrete. Therefore, the statements (a) and (b) imply (c). To prove the converse, we only need to recall that if $X$ is uniformly discrete, then $$\hbox{Lip}(X,A)=\hbox{lip}(X,A)=\ell^\infty(X, A)=C(\beta X, A).$$
Since $\beta X$ is compact and $A$ is amenable, $C(\beta X, A)$ is amenable; see \cite{tik}. Hence $\frak{F}(X, A)$ is amenable. 
%We have $lip^\circ d(X, A)=lip(\beta X, A)$. On the other hand, $Lip(\beta X, A)=lip(\beta X, A)^{**}$ is amenable. So $lip(\beta X, A)$ is amenable. So $lip^\circ d(X, A)$ is amenable.The other statements are proved similarly.
$\hfill\square$

\begin{theorem}  Let $(X, d)$ be a locally compact metric space and $A$ be a Banach algebra with $\Delta(A)\neq\emptyset$.  Then $\emph{lip}^\circ(X,A)$ is amenable if and only if $X$ is uniformly discrete and $A$ is amenable.
\end{theorem}
{\it Proof.} Let  $X$ be uniformly discrete and $A$ be amenable. Then
$$
\hbox{lip}^\circ(X,A)=c_0(X, A)=
c_0(X)\check{\otimes}A= \overline{c_0(X)\hat{\otimes}A}^{\|.\|_\epsilon}.
$$
Thus $\hbox{lip}^\circ(X,A)$ is amenable. The converse is proved similar to that of the proof of Theorem \ref{poi}. So we omit it.$\hfill\square$\\

Now, we discussed the amenability the second dual of Lipschitz algebras. For study of the second dual of Lipschitz algebras see \cite{bcd, biyabani1}.

\begin{theorem}Let $X$ be a metric space, $A$ be a Banach algebra with $\Delta(A)\neq\emptyset$ and $A^*$ be separable. Then the following assertions are equivalent.

\emph{(a)} $\emph{Lip}(X, A)^{**}$ \emph{(}respectively, $\emph{lip}(X, A)^{**}$\emph{)} is amenable.

\emph{(b)} $\emph{Lip}(X, A)$ \emph{(}respectively, $\emph{lip}(X, A)$\emph{)} is amenable and $A^{**}$ is amenable

\emph{(c)} $X$ is uniformly discrete and $A^{**}$ is amenable.\\
In this case $A^{**}$ and $\emph{Lip}(X, A)^{**}$ \emph{(}respectively, $\emph{lip}(X, A)^{**}$\emph{)} are unital.
\end{theorem}
{\it Proof.} Let $\hbox{Lip}(X, A)^{**}$ be amenable. It follows from Theorem 2.3 in \cite{gstud} that $\hbox{Lip}(X, A)$ is amenable and so $X$ is uniformly discrete. On the other hand, there exists a continuous epimorphism from $\hbox{Lip}(X, A)^{**}$ onto $A^{**}$. Hence amenability of $\hbox{Lip}(X, A)^{**}$ implies amenability of $A^{**}$. Thus (a)$\Rightarrow$(b)$\Rightarrow$(c).
If $X$  is uniformly discrete and $A^{**}$ is amenable, then $C(\beta X)^{**}\hat{\otimes}A^{**}$ is amenable. Since $A^*$ is separable,
$$
\hbox{Lip}(X, A)^*=C(\beta X, A)^*=M(\beta X, A^*)=M(\beta X)\hat{\otimes}A^*.
$$
Thus
 $$\hbox{Lip}(X,A)^{**}=\overline{C(\beta X)^{**}\hat{\otimes}A^{**}}^{\|.\|_\epsilon}$$
is amenable. That is, (c) $\Rightarrow$(a). The other case is proved similarly.

Finally, note that if the second dual of a Banach algebra is amenable, then it  is unital; see Lemma 1.1 in \cite{glw}.$\hfill\square$

\begin{corollary}\label{0124} Let $X$ be a metric space.  Then the following assertions are equivalent.

\emph{(a)} $(\emph{Lip}X)^{**}$ is amenable.

\emph{(b)} $\emph{Lip}X$ is amenable.

\emph{(c)} $(\emph{lip}X)^{**}$ is amenable.

\emph{(d)} $\emph{lip}X$ is amenable.

\emph{(e)} $X$ is uniformly discrete.
\end{corollary}

\begin{theorem}Let $X$ be a locally compact metric space and $A$ be a Banach algebra with $\Delta(A)\neq\emptyset$. Then $\emph{lip}^\circ(X, A)^{**}$ is amenable if and only if $X$ is finite and $A^{**}$ is amenable.
\end{theorem}
{\it Proof.} Let $\hbox{lip}^\circ(X, A)^{**}$ is amenable. It follows from Theorem 2.3 in \cite{gstud} and Lemma 1.1 in \cite{glw} that $\hbox{lip}^\circ(X, A)$ is a unital amenable Banach algebra. 
This together with Theorem \ref{poi} shows that $G$ is finite. For every $x\in X$, the mapping $f\mapsto x\otimes f$ is a continuous epimorphism from $\hbox{lip}^\circ(X, A)$ onto $A$. So from amenability of $\hbox{lip}^\circ(X, A)^{**}$ we see that $A^{**}$ is amenable. The converse follows from the fact that if $X$ is finite, then $\hbox{lip}^\circ(X, A)=\sum_{i=1}^n\Bbb{C}\hat{\otimes}A^{**}$ for some $n\in\Bbb{N}$. $\hfill\square$\\

The following result is needed in the sequel.

\begin{proposition}\label{un12} Let $A_1$ and $A_2$ be Banach algebras. Then $A_1$ and $A_2$ are amenable if and only if $A_1\hat{\otimes}A_2$ is amenable.
\end{proposition}
{\it Proof.} Let $\pi_i: A_1\times A_2\rightarrow A_i$ be canonical projective on $A_i$. Then $\pi_i$ is bilinear and surjective. Hence there exists a continuous epimorphism from $A_1\hat{\otimes}A_2$ onto $A_i$. Now, if $A_1\hat{\otimes}A_2$ is amenable, then $A_1$ and $A_2$ are amenable. The converse follows from the fact that approximate diagonals for $A_1$ and $A_2$ construct an approximate diagonal for $A_1\hat{\otimes}A_2$.$\hfill\square$\\

Chameh and the second author \cite{cr} proved that the injective tensor product $\hbox{Lip}(X)\check{\otimes}A$ can be embedded isometrically into $(\hbox{Lip}(X, A), \||.\||)$ and so it is a closed subalgebra of $(\hbox{Lip}(X, A), \||.\||)$, where $\||f\||=\max\{\|f\|_{\infty, A}, \|f\|_{d, A}\}$. Clearly, $\|.\|_{d, A}$ and $\||.\||$ are equivalent.

\begin{theorem}\label{asab}  Let $(X, d)$ be a metric space and $A$ be a Banach algebra with $\Delta(A)\neq\emptyset$. Then the following assertions are equivalent.

\emph{(a)} $\emph{Lip}(X)^{**}\check{\otimes}A^{**}$ \emph{(}respectively, $\emph{lip}(X)^{**}\check{\otimes}A^{**}$\emph{)} is amenable.

\emph{(b)} $\emph{Lip}(X)^{**}\hat{\otimes}A^{**}$ \emph{(}respectively, $\emph{lip}(X)^{**}\hat{\otimes}A^{**}$\emph{)} is amenable.

\emph{(c)} $X$ is uniformly discrete and $A^{**}$ is amenable.
\end{theorem}
{\it Proof.} Let $\varphi\in\Delta(A)$. Then the mapping $F\otimes m\mapsto\varphi^{**}(m) F$ is a continuous epimorphism from 
 $\hbox{Lip}(X)^{**}\check{\otimes}A^{**}$ onto $\hbox{Lip}(X)^{**}$. So if (a) holds, then $\hbox{Lip}(X)^{**}$ is amenable and hence $X$ is uniformly discrete. Similarly, $A^{**}$ is amenable. Thus (a)$\Rightarrow$(c). Since
$$
\hbox{Lip}(X)^{**}\check{\otimes}A^{**}= \overline{\hbox{Lip}(X)^{**}\hat{\otimes}A^{**}}^{\|.\|_\epsilon},
$$
we have (b)$\Rightarrow$(a). By Corollary \ref{0124} and proposition \ref{un12}, (c)$\Rightarrow$(b).
$\hfill\square$\\

By the same argument as used in the proof of Theorem \ref{asab} shows that $\hbox{lip}^\circ(X)^{**}\check{\otimes}A^{**}$  is amenable
if and only if $\hbox{lip}^\circ(X)^{**}\hat{\otimes}A^{**}$  is amenable; or equivalently,  $X$ is finite and $A^{**}$ is amenable, where $(X, d)$ is a locally compact metric space and $A$ is a Banach algebra with $\Delta(A)\neq\emptyset$.

\section{\normalsize\bf Weak amenability of $\hbox{Lip}(X,A)$}

Now, we prove some results concerning  point amenability vector-valued Lipschitz algebras.

\begin{theorem}\label{efi1}  Let $(X, d)$ be a metric space and $A$ be a Banach algebra. Then the following statements hold.

\emph{(i)} If $\emph{Lip}(X,A)$ is point amenable, then  $X$ is discrete and $A$ is point amenable.

\emph{(ii)} If $\emph{Lip}(X,A)$ is weakly amenable or cyclically weakly amenable, then  $X$ is discrete and $A$ is point amenable.
\end{theorem}
{\it Proof.} Suppose that $X$ is non-discrete. Then there exists $x_0\in X$ such that for any $n\in\Bbb{N}$, we have $B_d(x_0,1/n)\neq\{x_0\}$. Hence there exists a sequence $(x_n)$ such that $$0<d(x_n,x_0)<1/n.$$ It follows that  $x_n\rightarrow x_0$. Fix $\varphi\in\Delta(A)$ and for every $n\in\Bbb{N}$ define the linear functional $Q_n$ on $\hbox{Lip}(X, A)$  by $$Q_n(f)= \frac{\varphi(f(x_n))-\varphi(f(x_0))}{d( x_n,x_0)}.$$ Then
$$
|Q_n(f)|\leq
\frac{\|f(x_n)-f(x_0)\|}{d( x_n, x_0)}\leq\|f\|_{d, A}.
$$
It follows that $(Q_n)$ is a sequence in the unit ball of $\hbox{Lip}(X, A)^*$.  So $Q_n\rightarrow Q$ in the weak$^*$ topology of $\hbox{Lip}(X, A)^*$ for some $Q\in\hbox{Lip}(X, A)^*$. Thus for every $f, g\in\hbox{Lip}(X, A)$, we have
$$
Q_n(fg)=\varphi(f( x_n))Q_n(g)+\varphi(g( x))Q_n(f) .
$$
This shows that
$$
Q( fg)=x_0\otimes\varphi (f)Q(g)+x_0\otimes\varphi(g)Q( f).
$$
Since $\varphi$ is non-zero, there exists $a\in A$ such that $\varphi(a)=1$. Define the vector-valued function $h\in\hbox{Lip}(X, A)$ by
 $$
h(y)=\min\{1, d(x_0, y)\}a.
$$
Then $h(x_0)=0$ and $h(x_n))=d(x_n,x_0)a$ for all $n\in\Bbb{N}$. Thus for every $n\in\Bbb{N}$, we have $Q_n(h)=1$ and so $Q(h)=1$. Therefore, $Q$ is a non-zero continuous point derivation of $\hbox{Lip}(X, A)$. That is, $\hbox{Lip}(X, A)$ is not 0-point amenable. So (i) holds. The statement (ii) follows from (i) and Theorem 4.1 in \cite{mr4}.$\hfill\square$\\

From Theorem 4.6 in \cite{mr4} and Theorem \ref{efi1} we see that the following result holds.

\begin{corollary} Let $(X, d)$ be a metric space and $A$ be a unital, commutative Banach algebra. Then the following assertions are equivalent.

\emph{(a)} $\emph{Lip}(X,A)$ is weakly amenable.

\emph{(b)} $\emph{Lip}(X,A)$ is cyclically weakly amenable.

\emph{(c)} $\emph{Lip}(X,A)$ is point amenable.

\emph{(d)}  Every maximal ideal of $\emph{Lip}(X,A)$ is essential.\\
In this case, $X$ is discrete and $A$ is point amenable.
\end{corollary}

The next example shows that Theorem \ref{efi1} does not remain valid for $\hbox{lip}(X, A)$ instead of $\hbox{Lip}(X, A)$.

\begin{example}{\rm (i) Let $d$ be a metric function on  $\Bbb{T}$ be a metric space and $\alpha\in(0, 1/2)$. Then $(\Bbb{T}, d^\alpha)$ is also a metric space and so the corresponding the Lipschitz algebra $\hbox{lip}_\alpha\Bbb{T}$ is weakly amenable, however, $\Bbb{T}$ is not discrete.

(ii) Let $X=\Bbb{N}$ and $A$ be the unital, commutative Banach algebra $\hbox{lip}_{2/3}(\Bbb{T})$. Then $A$ is not weakly amenable and by Theorem 2.1 (iv) in \cite{mr5}, the Banach algebra $\hbox{Lip}(X, A)$ is not weakly amenable, however, $X$ is discrete. So the converse of Theorem \ref{efi1} does not hold, in general.}
\end{example}

For a metric space $X$ and Banach algebra $A$, let $$X\otimes\Delta(A)=\{x\otimes\varphi: x\in X, \varphi\in\Delta(A)\}$$
The following result is an analogue of Proposition 9.2  in \cite{s2}.

\begin{proposition}\label{efi2}  Let $(X, d)$ be a discrete metric space and $A$ be a point amenable Banach algebra. Then every continuous point derivation at $X\otimes\Delta(A)$ of $\emph{Lip}(X,A)$ is zero.
\end{proposition}
{\it Proof.} Let $\varphi\in\Delta(A)$ and $x_0\in X$. Assume that $d$ is a continuous point derivation at $x_0\otimes\varphi$ of $A$.
Define the linear functional $\tilde{d}_{x_0}$ on $A$ by $$\tilde{d}_{x_0}(a)=d(\chi_{x_0}\otimes a).$$ Then $\tilde{d}_{x_0}$ is a point derivation at $\varphi$ of $A$. Since $A$ is point amenable, $\tilde{d}_{x_0}=0$ and so $d(\chi_{x_0}\otimes a)=0$. That is, $d=0$ on $\langle\{\chi_{x_0}\otimes a: a\in A\}\rangle$. Let $x\neq x_0$. Choose $a_0, a\in A$ with $\varphi(a_0)=1$ and $\varphi(a)\neq 0$. Set $Q=x_0\otimes\varphi$, $$f=\chi_{x}\otimes a\quad\hbox{and}\quad g=\chi_{x_0}\otimes a_0.$$
Then $fg=0$. Since $Q(f)=0$ and $Q(g)=1$, it follows that
$$
d(f)=Q(f)d(g)+Q(g)d(f)=d(fg)=0.
$$
Thus $d(\chi_x\otimes a)=0$ for all $x\in X$ and $a\in A$.
But, $\hbox{Lip}(X, A)$ is a subspace of $\ell^\infty(X, A)$. Hence $d=0$. In fact, if $f\in\hbox{Lip}(X, A)$, then  $f\in\ell^\infty(X)\check{\otimes}A$. Thus there exists a sequences $(f_n)_{n\in\Bbb{N}}$ in $\ell^\infty(X)\otimes A$ such that $f_n\rightarrow f$. We have $d(f_n)=0$ for all $n\in\Bbb{N}$. This shows that $d=0$.$\hfill\square$\\

%
% Since $X$ is discrete,
%$$
%\alpha:=\inf_{x\in X, x\neq {x_0}}d(x, {x_0})>0.
%$$
%Choose $a\in A$ with $\varphi(a)\neq 0$. Then
%\begin{eqnarray*}
%\sup_{x,y\in X, x\neq y}\frac{\|a\chi_{x_0}( x)-a\chi_{x_0}(y)\|}{d(x,y)}=
%\|a\|/\alpha<\infty,
%\end{eqnarray*}
%where $\chi_{x_0}$ denotes the characteristic function of $\{x_0\}$ on $X$. This shows that $g:=\chi_{x_0}\otimes a\in\hbox{Lip}(X, A)$. If $f\in \hbox{Lip}(X, A)$, then for every $y\in X$, we have
%\begin{eqnarray*}
%(gf)(y)&=&g(y)f(y)\\
%&=&a\chi_{x_0}(y)f(y)\\
%&=&a\chi_{x_0}(y)f(x_0)\\
%&=&f({x_0})g(y).
%\end{eqnarray*}
%Hence
%$$
%gf=f(x_0)g.
%$$
%Now, let $d$ be a continuous derivation at $Q:=\delta_{x_0}\otimes\varphi$ of $\hbox{Lip}(X, A)$. Then
%\begin{eqnarray}\label{chi}
%f(x_0)d(g)&=&d(gf)\nonumber\\
%&=&Q(g)d(f)+Q(f)d(g)\nonumber\\
%&=&\delta_{x_0}(\chi_{x_0})\varphi(a)d(f)+Q(f)d(g)\\
%&=&\varphi(a)d(f)+Q(f)d(g).\nonumber
%\end{eqnarray}
%We define the linear functionals $q$ and $\tilde{d}$ in $A^*$ by $$
%q(a)=Q(\chi_{x_0}\otimes a)\quad\hbox{and}\quad\tilde{d}(a)=d(\chi_{x_0}\otimes a).$$  Then $q$ and $\tilde{d}$ are well-defined and for every $a, b\in A$, we have
%$$
%\tilde{d}(ab)= q(a)\tilde{d}(b)+q(b)\tilde{d}( a).
%$$
%Since $A$ is point amenable and $\tilde{d}$ is a continuous point derivation at $q$ of $A$, it follows that $\tilde{d}=0$. So $d(g)=0$.
%This together with (\ref{chi}) implies that $d(f)=0$. That is, $d=0$ on $\hbox{Lip}(X, A)$. $\hfill\square$\\

As an immediate consequence of Theorem \ref{efi1} and Proposition \ref{efi2} we have the following result.

\begin{corollary} Let $(X, d)$ be a metric space and $A$ be a Banach algebra such that $\Delta(\emph{Lip}(X, A))=X\otimes\Delta(A)$. Then $\emph{Lip}(X,A)$ is point amenable if and only if  $X$ is discrete and $A$ is point amenable.
\end{corollary}

\begin{corollary} Let $(X, d)$ be a compact metric space and $A$ be a  Banach algebra. Then $\emph{Lip}(X,A)$ is point amenable if and only if  $X$ is finite and $A$ is point amenable.
\end{corollary}

\section{\normalsize\bf Weak amenability of $\hbox{Lip}X$}

%Let $X$ and $Y$ be locally compact spaces. A bounded linear functional on $C_0(X)\hat{\otimes}C_0(Y)$ is called a \emph{bimeasure on} $X\times Y$.
For metric space $(X, d)$ and $\alpha\in (0, 1]$, Sherbert proved that if $x$ is a cluster
point of $X$, then the zero map is the only point derivation at $\delta_x$ of $\hbox{Lip}_\alpha X$, where $\delta_x(f)=f(x)$ for all $f\in \hbox{Lip}_\alpha X$; see Proposition 9.2 in \cite{s2}. This implies that if $\hbox{Lip}_\alpha X$ is weakly amenable, cyclically weakly amenable or point amenable, then $X$ is discrete. We show that the converse holds. Before, we give the proof of it, let us recall that  a  bounded linear functional on $C_0(X)\hat{\otimes}C_0(Y)$ is called a \emph{bimeasure on} $X\times Y$, where $X$ and $Y$ be locally compact spaces. The space of all bimeasures on $X$ is denoted by $BM(X, Y)$. If $F\in BM(X, Y)$, then there exist positive constant $k$ and regular Borel probability measures $\mu$ on $X$ and $\nu$ on $Y$ such that 
$$
|\langle F, f\otimes g\rangle|\leq k\|f\|_2\|g\|_2
$$
for all $f\in L^2(X, \mu)$ and $g\in L^2(Y, \nu)$. 
Also, $F$ can be extended to  a bounded linear functional on $$L^2(X, \mu)\otimes_2  L^2(Y, \nu)=L^2(X\times Y, \mu\times\nu).$$
Note that
$$
\|\sum_{i=1}^n f_i\otimes g_i\|_2=(\int_X\int_Y|\sum_{i=1}^nf_i(x)g_i(y)|^2d\mu(x)d\nu(y))^{1/2},
$$
where $f_i\in L^2(X, \mu)$, $g_i\in L^2(Y, \nu)$ and $n\in\Bbb{N}$. 
Thus $$F\in L^2(X\times Y, \mu\times\nu)^*=L^2(X\times Y, \mu\times\nu).$$
So there exists a function $w\in L^2(X\times Y, \mu\times\nu)$ such that $$\langle F, f\otimes g\rangle=\langle f\otimes g, w\rangle_2.$$
Therefore,
$$
\langle F, f\otimes g\rangle=\int_X\int_Yf(x)g(y)w(x, y) d\mu(x)d\nu(y).
$$
One can prove that $B$ is finite dimensional and $A$ be a Banach algebra, then $A\check{\otimes}B=A\hat{\otimes}B$; see for example \cite{r1}. So  if $Y$ is finite dimensional, then
\begin{eqnarray*}
BM(X, Y)&=&(C_0(X)\check{\otimes}C_0(Y
))^*\\
&=&C_0(X\times Y)^*\\
&=&M(X\times Y)\\
&=&M(X)\hat{\otimes}M(Y)\\
&=&M(X)\hat{\otimes}C_0(Y)^*\\
&=&M(X)\check{\otimes}C_0(Y)\\
&=&
M(X)\hat{\otimes}C_0(Y)=M(X, C_0(Y)).
\end{eqnarray*}

\begin{theorem}\label{an} Let $(X, d)$ be a metric space. Then the following assertions are equivalent.

\emph{(a)} $\emph{Lip}X$ is weakly amenable.

\emph{(b)} $\emph{Lip}X$ is cyclically weakly amenable.

\emph{(c)} $\emph{Lip}X$ is point amenable.

\emph{(d)} $X$ is discrete.
\end{theorem}
{\it Proof.} Let $X$ be discrete. Then $X$ is a locally compact Hausdorff space. By the Stone-Weierstrass theorem, $\hbox{Lip} X$ is a dense subspace of $\ell^\infty(X)$. So the mapping $$Q:=\iota\otimes\iota: \hbox{Lip} X\hat{\otimes}\hbox{Lip} X\rightarrow \ell^\infty(X)\hat{\otimes}\ell^\infty(X)$$ is a continuous homomorphism with dense range, where $\iota: \hbox{Lip} X\rightarrow \ell^\infty(X)$ is the inclusion map. Assume now that $D: \hbox{Lip} X\rightarrow (\hbox{Lip} X)^*$ is a continuous derivation. Then there exists
$T\in (\hbox{Lip} X\hat{\otimes}\hbox{Lip} X)^*$ such that
$$
\langle D(f), g\rangle=\langle T, f\otimes g\rangle
$$
for all $f, g\in \hbox{Lip} X$; see for example  Proposition 13 VI in
\cite{bd}.
Since $Q^*$ is surjective, $Q^*(F)=T$ for some $F\in(\ell^\infty(X)\hat{\otimes}\ell^\infty(X))^*$. But $\ell^\infty(X)$ is a $C^*-$algebra. Consequently, there exist $\mu, \nu\in M_b(\beta X)$, constant $M>0$ and bounded bilinear map $\bar{F}: L^2(\beta X, \mu)\otimes_2 L^2(\beta X, \nu)\rightarrow\Bbb{C}$ such that $$\bar{F}|_{\ell^\infty(X)\hat{\otimes}\ell^\infty(X)}=F
$$
and $|\bar{F}(f, g)|\leq M\|f\|_2\|g\|_2$ for all $f\in L^2(\beta X, \mu)$ and $g\in L^2(\beta X, \nu)$. So 
$$
|T(f\otimes g)|=|\langle D(f), g\rangle|\leq M\|f\|_2\|g\|_2$$ for all $f, g\in \ell^2(X)$; see Corollary 1.3 in \cite{gs}. This shows that
$$T\in (\ell^2(X)\otimes_2\ell^2(X))^*=\ell^2(X\times X)^*.
$$
Thus there exists $k\in \ell^2(X\times X)$ such that for every $f, g\in \ell^2(X)$
$$
\langle D(f), g\rangle=\int_X\int_Y f(x) g(y)\; dw(x, y),
$$
where $$dw(x, y)=k(x, y)d\mu(x)d\nu(y).$$ Since $D$ is a derivation, $D(1)=0$ and so
$$
\int_X\int_Y dw(x, y)=\langle D(1), 1\rangle=0.
$$
Therefore, $w=0$. This implies that $D=0$. That is, (d) implies (a).
From Proposition 9.2 in \cite{s2} we see that (a) implies (d). Since $\hbox{Lip}X$ is a unital commutative Banach algebra, by Theorem 4.6 in \cite{mr4}, the statements (a), (b) and (c) are equivalent.$\hfill\square$

\begin{remark} {\rm Let $X$ be $\Bbb{R}$ with the usual Euclidean metric and $\alpha>1$. Then 
$$
\hbox{Lip}_\alpha X=\hbox{Con}(X)=\Bbb{C},
$$
where $\hbox{Con}(X)$ is the space
consisting of all constant functions on $X$. This shows that $\hbox{Lip}_\alpha X$ is weakly amenable, however, $X$ is not discrete. So Theorem \ref{an} dose not remain valid for $\hbox{Lip}_\alpha$ when $\alpha>1$.}
\end{remark}

\begin{remark} {\rm We give another proof for the implication (d)$\Rightarrow$(c) of Theorem \ref{an}. Let $X$ be discrete. Then $\hbox{Lip}X$ is unital and every maximal ideal $M$ in $\hbox{Lip}X$ has the form
$$
M=\{f\in\hbox{Lip}X: \hat{f}(z_0)=0\}
$$
for some $z_0\in l_X:= \overline{\{\delta_x: x\in X\}}^{w^*}$. For every $n\in\Bbb{N}$ and $x\in X$, let
$$
e_n(x)=\min\{n\|\delta_x-\delta_{z_0}\|, 1\}.
$$
Then for every $x\in X$ there exists $n_0\in\Bbb{N}$ such that $e_n(x)=1$ for all $n\geq n_0$. Indeed, if $f=\chi_{z_0}$, then $f\in \hbox{Lip}X$ and $\|f\|_{d, \Bbb{C}}=1+1/a$, where $$a=d(z_0, l_X-\{z_0\})>0.$$ So if $x\neq z_0,$ then
$$
\|\delta_x-\delta_{z_0}\|\geq\frac{|\delta_{z_0}(x)-\delta_{z_0}(z_0)|a}{1+a}=\frac{a}{1+a}.
$$
Thus $e_n(x)=1$ for all $n>\frac{1+a}{a}$. So $fe_n=f$ and hence $$\|fe_n-f\|_{d, \Bbb{C}}=0$$ for all $n>\frac{1+a}{a}$.
This shows that  $(e_n)$ is a bounded approximate identity for $M$. So $M=M^2$. This implies that  $\hbox{Lip}X$ is point amenable; see Theorem 4.6 in \cite{mr4}. }
\end{remark}

\begin{corollary} Let $(X, d)$ be a compact metric space. Then the following assertions are equivalent.

\emph{(a)} $\emph{Lip}X$ is weakly amenable.

\emph{(b)} $\emph{Lip}X$ is cyclically weakly amenable.

\emph{(c)} $\emph{Lip}X$ is point amenable.

\emph{(d)} $X$ is finite.
\end{corollary}

As another consequence of Theorem \ref{an}, we have the following result due to Sherbert \cite{s2}.

\begin{corollary} Let $(X, d)$ be an infinite compact metric space and $\alpha\in (0, 1]$. Then the following statements hold.

\emph{(i)}  $\emph{Lip}_\alpha X$ is not point amenable.

\emph{(ii)} $\emph{Lip}_\alpha X$ is not weakly amenable.

\emph{(iii)} $\emph{Lip}_\alpha X$ is not cyclically weakly amenable.
\end{corollary}

\begin{corollary} Let $(X, d)$ be a metric space.  If $d$ is bounded below, then $\emph{Lip}X$ is weakly amenable.
\end{corollary}
{\it Proof.} Let $d$ be bounded below. Then there exists $M>0$ such that $d(x, y)>M$ for all $x, y\in X$ and $x\neq y$. So $X$ is discrete.  By Theorem \ref{an}, $\hbox{Lip}X$ is weakly amenable.$\hfill\square$

\section{\normalsize\bf Weak amenability of vector-valued Banach algebras}

First, we investigate cohomological properties of Banach algebra $C_0(X, A)$.

\begin{theorem}\label{havig2} Let $X$ be a locally compact Hausdorff space and $A$ be a Banach algebra. Then the following statements hold.

\emph{(i)} $C_0(X, A)$ is amenable if and only if $A$ is amenable. 

\emph{(ii)} $C_0(X, A)$ is 0-point amenable if and only if $A$ is 0-point amenable.

\emph{(iii)} If $A$ is commutative, then $C_0(X, A)$ is weakly amenable if and only if $A$ is weakly amenable. 
\end{theorem}
{\it Proof.} Let $C_0(X, A)$ be (respectively, 0-point) amenable. Since continuous epimorphisms are (respectively, 0-point) amenability preserving, $A$ is (respectively, 0-point); see Proposition 2.3.1 in \cite{r1} and Theorem 2.1 in \cite{mr5}. For the converse, note that from Proposition 1.5.6 in \cite{k} and Lemma \ref{havig1} we see that
$$
C_0(X,A)=
C_0(X)\check{\otimes}A= \overline{C_0(X)\hat{\otimes}A}^{\|.\|_\epsilon}.
$$
So if $A$ is (respectively, 0-point) amenable, then Theorem 4.1 in \cite{mr5} and Proposition \ref{un12} show that $C_0(X)\hat{\otimes}A$ is (respectively, 0-point) amenable. Therefore, $C_0(X,A)$ is (respectively, 0-point) amenable. Hence (i) and (ii) are proved. 

Now, let $A$ be commutative. If $C_0(X, A)$ be weakly amenable, then $C_0(X,A)$ is commutative and cyclically weakly amenable. In view of Theorem 2.1 in \cite{mr5}, $A$ is weakly amenable. To prove the converse, let us remark from \cite{gro} that  $C_0(X)\hat{\otimes}A$ is weakly amenable. By Corollary 2.3 in \cite{mr5}, the Banach algebra $C_0(X,A)$ is weakly amenable.$\hfill\square$\\

As a consequence of Theorem 4.6 in \cite{mr4} and Theorem \ref{havig2}
we have the next result.

\begin{corollary}\label{anar} Let $(X, d)$ be a metric space and let $A$ be a  unital commutative Banach algebra. Then the following assertions are equivalent.

\emph{(a)} $C_0(X, A)$ is point amenable.

\emph{(b)} $C_0(X, A)$ is cyclically weakly amenable.

\emph{(c)} $C_0(X, A)$ is weakly amenable.

\emph{(d)}  $A$ is weakly amenable.
\end{corollary}

\begin{remark}{\rm Let $X$ be a locally compact Hausdorff spaces and let $G$ be a locally compact  group. It is well-known from \cite{j} that $L^1(G)$ is weakly amenable and so it is 0-point amenable. Also, $M(G)$ is point amenable if and only if $M(G)$ is weakly amenable; or equivalently, $G$ is discrete. From these facts and Theorem \ref{havig2} we infer that 
$C_0(X, L^1(G))$ is 0-point amenable, however, 0-point amenability of 
$C_0(X, M(G))$ is equivalent to discreteness of $G$. In the case where $G$ is abelian, Corollary \ref{anar} implies that $C_0(X, M(G))$ is weakly amenable if and only if $G$ is discrete.}
\end{remark}

For a locally compact group $G$ and Banach algebra $A$, let $L^1(G, A)$ be the Banach algebra of $A-$valued integrable
functions on $G$. Note that $L^1(G, A)$ is isometrically isomorphic to $L^1(G)\hat{\otimes}A$; see for example Proposition 1.5.4 in [19].

\begin{theorem}\label{paran} Let $G$ be a locally compact group and $A$ be a Banach algebra. Then the following statements hold.

\emph{(i)} $L^1(G, A)$ is amenable if and only if $G$ and $A$ are amenable.

\emph{(ii)} $L^1(G, A)$ is point amenable if and only if $A$ is point amenable.

\emph{(iii)} Assume that $G$ is abelian and $A$ is commutative. Then $L^1(G, A)$ is weakly amenable if and only if $A$ is weakly amenable.
\end{theorem}
{\it Proof.}
% In view of Proposition 1.5.4 in \cite{k},
%\begin{eqnarray}\label{anar1}
%L^1(G, A)=L^1(G)\hat{\otimes}A,
%\end{equanrray}
Groenbaek \cite{gro} proved that the projective tensor product two  weakly amenable commutative Banach algebras is weakly amenable. Yazdanpanah \cite{y} showed that the converse of this result holds.
The authors \cite{mr5} established the Groenbaek's and Yazdanpanah's   result for 0-point amenability instead of weakly amenability. Theses facts together with Proposition 1.5.4 in \cite{k} and Proposition \ref{un12} prove the theorem.
$\hfill\square$

\begin{example}{\rm Let  $X$ be a locally compact Hausdorff space and $G$ be a  locally compact group. Let $\frak{A}$ be one of the Banach algebras   $L^1(G, L^1(G))$ or $L^1(G, C_0(X))$. Then  by Theorem \ref{paran}, $\frak{A}$ is point amenable. Also
if $G$ is abelian, then $\frak{A}$ is weakly amenable.}
\end{example}

\begin{corollary}\label{tala} Let $G_i$ be a locally compact abelian group and $A_i$ be a commutative Banach algebra, for  $i=1,2$. Then the following assertions are equivalent.

\emph{(a)}  $L^1(G_1, A_1)\hat{\otimes}L^1(G_2, A_2)$ is weakly amenable.

\emph{(b)} $L^1(G_1\times G_2, A_1\hat{\otimes}A_2)$  is weakly amenable.

\emph{(c)} $A_1$ and $A_2$ are weakly amenable.
\end{corollary}
{\it Proof.} In view of Proposition 1.5.4 in \cite{k}, 
$$
L^1(G_i, A_i)=L^1(G_i)\hat{\otimes} A_i
$$ 
for $i=1, 2$. It follows that
\begin{eqnarray*}
L^1(G_1,A_1)\hat{\otimes}L^1(G_2, A_2)&=&
L^1(G_1)\hat{\otimes}A_1\hat{\otimes}L^1(G_2)\hat{\otimes}A_2\\
&=&L^1(G_1\times G_2, A_1\hat{\otimes}A_2).
\end{eqnarray*}
This together with Theorem \ref{paran} proves the corollary.$\hfill\square$

\begin{corollary}  Let $G_i$ be a locally compact abelian group and $A_i$ be a commutative unital Banach algebra, for  $i=1,2$. Then the following assertions are equivalent.

\emph{(a)}  $L^1(G_1, A_1)\hat{\otimes}L^1(G_2, A_2)$ is  weakly amenable.

\emph{(b)}  $L^1(G_1, A_1)\hat{\otimes}L^1(G_2, A_2)$ is cyclically weakly amenable.

\emph{(c)}  $L^1(G_1, A_1)\hat{\otimes}L^1(G_2, A_2)$ is 0-point  amenable.

\emph{(d)} $A_1$ and $A_2$ are weakly amenable.
\end{corollary}

Let us recall that a Banach space $E$ has the {\it Radon-Nikodym property} if for every bounded subset $Y$ of $E$ and $\epsilon>0$ there exists $y\in Y$ such that $y\not\in\overline{\langle Y\setminus B_\epsilon\rangle}$, where $B_\epsilon=\{z\in E: \|x-z\|<\epsilon\}$. Note that reflexive spaces and the spaces $\ell^1(E)$ for a set $E$ have the Radon-Nikodym 
property. 
One can prove that if $A$ has Radon-Nikodym property, then $M(G)\hat{\otimes}A=M(G, A)$; see for example \cite{du}. 

\begin{theorem}\label{mina} Let $G$ be a locally compact group and $A$ be a Banach algebra with $\Delta(A)\neq\emptyset$. If $A$ has the Radon-Nikodym property, then the following assertions are equivalent.

\emph{(a)} $M(G, A)^{**}$ is 0-point amenable.

%\emph{(b)} $G$ is discrete and $M(G, A)^{**}$ is 0-point amenable.

\emph{(b)} $L^1(G, A)^{**}$ is 0-point amenable.

\emph{(c)} $G$ is finite and $A^{**}$ is 0-point amenable.
\end{theorem}
{\it Proof.} Let $M(G, A)^{**}$ is 0-point amenable. Then $(M(G)\hat{\otimes}A)^{**}$ be 0-point amenable. Then 
 $M(G)\hat{\otimes}A$ is 0-point amenable. By Theorem 4.1 in \cite{mr5}, the Banach algebra $M(G)$ is 0-point amenable. It follows from \cite{dgh}  that $G$ is discrete and so 
$$
(M(G)\hat{\otimes}A)^{**}=(\ell^1(G)\hat{\otimes}A)^{**}=L^1(G, A)^{**}
$$
is point amenable. Thus (a) implies (b). Since there exist continuous epimorphisms from $L^1(G, A)^{**}$ onto $L^1(G)^{**}$ and $A^{**}$, 
0-point amenability of $L^1(G, A)^{**}$ imply that $L^1(G)^{**}$ and $A^{**}$ are 0-point amenable. Thus $G$ is finite. That is, (b) implies (c).
To complete the proof, note that if $G$ is finite, then there is $n\in\Bbb{N}$ such that $(M(G)\hat{\otimes}A)^{**}$ and $\Bbb{C}^n\hat{\otimes} A^{**}$ are isomorphic. Thus (d) implies (a).$\hfill\square$

\begin{corollary}\label{ss} Let $G$ be a locally compact group and $A$ be a Banach algebra with $\Delta(A)$ is a non-empty set. If $A$ has the Radon-Nikodym property, then the following assertions are equivalent.

\emph{(a)} $M(G, A)^{**}$ is weakly amenable.

%\emph{(b)} $G$ is discrete and $M(G, A)^{**}$ is  cyclically \emph{(}weakly\emph{)} amenable.

\emph{(b)} $L^1(G, A)^{**}$ is  weakly amenable.

\emph{(c)} $G$ is finite and $A^{**}$ is  weakly amenable.
\end{corollary}
{\it Proof.} If $M(G, A)^{**}$ is weakly amenable, then $(M(G)\hat{\otimes}A)^{**}$ is 0-point amenable. By Theorem \ref{mina}, $G$ is finite and so $(M(G)\hat{\otimes}A)^{**}=L^1(G, A)^{**}.$ Thus (a) implies (b). If $L^1(G, A)^{**}$ is  weakly amenable, then by Theorem \ref{mina}, $G$ is finite. So there exists $n\in\Bbb{N}$ such that $$L^1(G, A)^{**}={\Bbb C}^n\hat{\otimes} A^{**}.$$
Thus ${\Bbb C}^n\hat{\otimes} A^{**}$ is weakly amenable. This together with Theorem 4.3 in \cite{mr5} shows that $A^{**}$ is weakly amenable. That is, (b) implies (c). Clearly, (d) implies (a).$\hfill\square$\\

%\begin{theorem} Let $A$ be a commutative regular Banach algebra. If $A^{**}$ is not weakly amenable, then $\ell^\infty(X, A)$ is not weakly amenable for some a set $X$.
%\end{theorem}
%{\it Proof.} It follows from \cite{.} that there exists a continuous epimorphism $\phi:  \ell^\infty(X, A)\rightarrow A^{**}$ for some a set $X$. So if $A^{**}$ is not weakly amenable, then by Theorem \ref{.} the Banach algebra $\ell^\infty(X, A)$ is not weakly amenable.$\hfill\square$\\

%\begin{theorem}  Let $A_i$ be a Banach algebra and $\varphi_i\in\Delta(A_i)$, for $i=1, 2$. Then $d$ is a continuous point derivation at $\varphi_1\otimes\varphi_2$ of $A_1\hat{\otimes} A_2$ if and only if there exist continuous point derivations $d_1$ at $\varphi_1$ of $A_1$ and $d_2$ at $\varphi_2$ of $A_2$ such that $d=d_1\otimes\varphi_2+d_2\otimes\varphi_1$.
%\end{theorem}

 Essmaili et al. \cite{esr} raised the question of whether there exists a Banach algebra $A$ such that $A^{**}$ is biprojective, but $A$ is not biprojective. The next result answers to the question.  

\begin{theorem} Let $A$ be a Banach algebra. Then the following statements hold.

\emph{(i)} If $A^{**}$ is biprojective, then $A$ is biprojective.

\emph{(ii)} If $A^{**}$ is cyclically weakly amenable, then $A$ is cyclically weakly amenable.

\emph{(iii)} If $A$ is abelian and $A^{**}$ is weakly amenable, then $A$ is weakly amenable.

\emph{(iv)} If $A$ is abelian and $A^{**}$ is cyclically amenable, then $A$ is cyclically amenable.

\end{theorem}
{\it Proof.}  Let $A^{**}$ be biprojective. By Proposition 2.6 in \cite{mr6}, $$(A^\sharp)^{**}=A^{**}\oplus\Bbb{C}$$
is biprojective. Since $(A^\sharp)^{**}$ is unital, it is contractible and so $A^\sharp$  is contractiblle; see Theorem 2.8 in \cite{mr6}. In view of Theorem 2.6 in \cite{erm123}, $A$ is contractible. It follows that $A$ is biprojective. That is, (i) holds.

Let $D: A\rightarrow A^*$ be a continuous derivation. Then $D^{**}$ is a continuous derivation; see Lemma 2.2 in \cite{gstud}.  If $A^{**}$ is cyclically weakly amenable, then $D^{**}$ is cyclic. Since $A$ is a subspace of $A^{**}$, it follows that $D$ is cyclic.
So $A$ is cyclically weakly amenable. Thus (ii) is proved. The statement (iii) follows from (ii) and the fact that every commutative weakly amenable Banach algebra is cyclically weakly amenable. In order to prove (iv), assume that $\Phi, \Psi\in A^{**}$. Then there exist nets $(a_\alpha)$ and $(b_\beta)$ in $A$ such that $a_\alpha\rightarrow\Phi$ and $b_\beta\rightarrow \Psi$ in the weak$^*$ topology of $A^{**}$. Hence
\begin{eqnarray*}\label{bbb}
\lim_\alpha\lim_\beta\langle D(a_\alpha), b_\beta\rangle=\langle D^{**}(\Phi), \Psi\rangle.
\end{eqnarray*}
This shows that if $D$ is cyclic, then $D^{**}$ is cyclic and so $D^{**}$ is inner. This implies that $D^{**}=0$. Therefore, $D=0$.$\hfill\square$

\begin{theorem} Let $X$ be a locally compact group and $A$ be a Banach algebra with a bounded approximate identity. Then the following assertions are equivalent.

\emph{(a)} $C_0(X, A)$ is biflat.

\emph{(b)} $A$ is biflat.

\emph{(c)} $A$ is amenable.
\end{theorem}
{\it Proof.} Since $A$ has a bounded approximate identity, (b) and (c) are equivalent. On the other hand, $C_0(X, A)$ has a bounded approximate identity. So $C_0(X, A)$ is biflat if and only if $C_0(X, A)$ is amenable; or equivalently, $A$ is amenable.$\hfill\square$

\vspace{2mm}

 {\footnotesize
\noindent {\bf Mohammad Javad Mehdipour}\\
Department of Mathematics,\\ Shiraz University of Technology,\\
Shiraz
71555-313, Iran\\ e-mail: mehdipour@.ac.ir\\
{\bf Ali Rejali}\\
Department of Pure Mathematics,\\ Faculty of Mathematics and Statistics,\\ University of Isfahan,\\
Isfahan
81746-73441, Iran\\ e-mail: rejali@sci.ui.ac.ir\\
\end{document}